\documentclass[12pt]{amsart}
\usepackage{amsfonts}
\usepackage{amssymb}
\theoremstyle{definition}

\theoremstyle{remark}

\newcommand{\C}{\mathbb C}
\newcommand{\D}{\mathbb D}
\newcommand{\Sph}{\mathbb S}

\begin{document}

\centerline{\bf Comptes rendus de l'Academie bulgare des Sciences}

\centerline{\it Tome 35, No 2, 1982}

\vspace{0.3in}

\centerline{\large\bf NEARLY K\"AHLER MANIFOLDS OF CONSTANT }
\centerline{\large\bf ANTIHOLOMORPHIC SECTIONAL CURVATURE}

\vspace{0.2in}
\centerline{\bf G. T. Ganchev, O. T. Kassabov}

\vspace{0.1in}
\centerline{\it (Submited by Academician B. Petkanchin on October 27, 1981)}

\vspace{0.1in}

In this paper we consider nearly K\"ahler manifolds of pointwise constant antiholomorphic 
sectional curvature. An analogue of Schur's theorem is proved and a classification 
theorem for such manifolds is given.

Let $M$ be an almost Hermitian manifold with dim\,$M=2n$, a metric tensor $g$,
an almost complex structure $J$ and a curvature tensor $R$. A 2-plane $\alpha$ in the
tangential space $T_pM$, $p\in M$ is said to be holomorphic (antiholomorphic) if
$J\alpha=\alpha$ ($J\alpha\perp\alpha$). The sectional curvature  $K(\alpha,p)$ of
the 2-plane $\alpha$ in $T_pM$ with an orthonormal basis $\{x,y\}$ is given by the 
equality $K(\alpha,p)=R(x,y,y,x)$. The manifold $M$ is said to be  of pointwise
constant holomorphic (antiholomorphic) sectional curvature if the sectional curvature
$K(\alpha,p)$ of an arbitrary holomorphic (antiholomorphic) 2-plane $\alpha$ in $T_pM$
for every $p\in M$ does not depend on $\alpha$.
An almost Hermitian manifold $M$ is called an $RK$-manifold if
$R(x,y,z,u)=R(Jx,Jy,Jz,Ju)$ for all $x,y,z,u \in T_pM$, $p\in M$. The nearly 
K\"ahler manifolds are characterized by the equality $(\nabla_XJ)X=0$, where $X$
is an arbitrary vector field. Every nearly K\"ahler manifold is an $RK$-manifold \cite{G2}.
 
For K\"ahler manifolds the following analogue of the classical Schur theorem is well known:

{\bf Theorem.} Let $M$ be a connected K\"ahler manifold with dim\,$M\ge 4$. If $M$ is of
pointwise constant holomorphic sectional curvature $\mu(p)$, then $\mu$ is a constant.
Moreover, $M$ is of constant antiholomorphic sectional curvature $\mu/4$. 

Conversely, if the K\"ahler manifold $M$ is of pointwise constant antiholomorphic sectional
curvature $\mu/4$, then $M$ is of constant holomorphic sectional curvature $\mu$ \cite{ChO}.
Such a manifold is locally isometric to $\C^n$, $\C\mathbb P^n$ or $\C\D^n$.

For nearly K\"ahler manifolds the following statements are known:

{\bf Theorem,} \cite{NH}. If $M$ is a connected nearly K\"ahler manifold of pointwise
constant holomorphic sectional curvature $\mu(p)$ and dim\,$M\ge4$, then $\mu$ is a 
constant.

{\bf Theorem,} \cite{G1}. If $M$ is a nearly K\"ahler manifold of constant holomorphic
sectional curvature and dim\,$M\ge 4$, then $M$ is locally isometric to
$\C^n$, $\C\mathbb P^n$, $\C\D^n$ or $\Sph^6$.

We shall consider nearly K\"ahler manifolds of pointwise constant antiholomorphic 
sectional curvature. We need the following lemma:

{\bf Lemma,} \cite{GG}. Let $M$ be an almost Hermitian manifold and $T$ be a tensor of
type (0,4) in $T_pM$ satisfying the conditions:

1) $T(x,y,z,u)=-T(y,x,z,u)$;

2) $T(x,y,z,u)+T(y,z,x,u)+T(z,x,y,u)=0$;

3) $T(x,y,z,u)=-T(x,y,u,z)$;

4) $T(x,y,z,u)=T(Jx,Jy,Jz,Ju)$;

5) $T(x,y,y,x)=0$, where $\{x,y\}$ is a basis of an arbitrary holomorphic or
antiholomorphic 2-plane in $T_pM$.

Then $T=0$.

Let $R'(x,y,z,u)=R(x,y,Jz,Ju)$. We denote by $S$, $S'$ and $\tau$, $\tau'$ the
Ricci tensors and scalar curvatures with respect to the tensors $R$, $R'$ respectively.
If $\{ E_1,\hdots,E_{2n} \}$ is an arbitrary orthonormal frame field, we have
$$
	\sum_{i=1}^{2n}(\nabla_{E_i}R)(X,Y,Z,E_i)=(\nabla_XS)(Y,Z)-(\nabla_YS)(X,Z) \ ; \leqno (1)
$$ 

\vspace{-0.2in}
$$
	\sum_{i=1}^{2n}(\nabla_{E_i}S)(X,E_i)=\frac12X\tau     \leqno (2)
$$
for arbitrary vector fields $X,\,Y,\,Z$. The following identities are valid for a nearly 
K\"ahler manifold \cite{G2}:
$$
	\sum_{i=1}^{2n}(\nabla_{E_i}S')(X,E_i)=\frac12X\tau'  \ ;     \leqno (3)
$$

\vspace{-0.2in}
$$
	X(\tau-\tau')=0  \ ;    \leqno (4)
$$
$$
	2(\nabla_X(S-S'))(Y,Z)=(S-S')((\nabla_XJ)Y,JZ)+(S-S')(JY,(\nabla_XJ)Z) \ .  \leqno (5)
$$

Let the tensors $R_1$, $R_2$, $\psi$ be defined by
$$
	R_1(x,y,z,u)=g(y,z)g(x,u)-g(x,z)g(y,u)\ ;
$$
$$
	R_2(x,y,z,u)=g(Jy,z)g(Jx,u)-g(Jx,z)g(Jy,u)-2g(Jx,y)g(Jz,u) \ ;
$$
$$
	\begin{array}{r}
		\psi(x,y,z,u)=g(Jy,z)S(Jx,u)-g(Jx,z)S(Jy,u)-2g(Jx,y)S(Jz,u) \ \  \\
	             +g(Jx,u)S(Jy,z)-g(Jy,u)S(Jx,z)-2g(Jz,u)S(Jx,y) \ .
	\end{array}
$$

{\bf Proposition 1.} Let $M$ be an $RK$-manifold of pointwise constant antiholomorphic
sectional curvature $\nu$ and dim\,$M=2n\ge4$. Then the curvature tensor $R$ has the form
$$
	R=\frac16\psi +\nu R_1-\frac{2n-1}3\nu R_2    \leqno (6)
$$

\vspace{-0.1in}
\noindent and

\vspace{-0.2in}
$$
	3S'-(n+1)S=\frac1{2n}(3\tau'-(n+1)\tau)g \ ,   \leqno (7)
$$
$$
	\nu= \frac{(2n+1)\tau-3\tau'}{8n(n^2-1)}  \ .  \leqno (8)
$$

{\bf Proof.} Let $\{x,y\}$ be an orthonormal basis for an arbitrary antiholomorphic 2-plane
in $T_pM$. From the condition $R(x,y,y,x)=\nu$ immediately follows
$$
	S(x,x)-R(x,Jx,Jx,x)=2(n-1)\nu \ ,   \leqno (9)
$$
where $x$ is an arbitrary unit vector. Now let $T=R-(1/6)\psi-\nu R_1+((2n-1)/3)\nu R_2$.
From the given condition, (9) and the Lemma it follows (6). By direct computation we find
(7) and (8).

{\bf Theorem 2.} Let $M$ be a connected nearly K\"ahler manifold of pointwise constant
antiholomorphic sectional curvature $\nu$ and dim\,$M=2n>4$. Then $\nu$ is a constant.

{\bf Proof.} From (7), (3) and (2) it follows that $X(3\tau'-(n+1)\tau)=0$ for an arbitrary
vector field $X$. Let $n>2$. Taking into account (4), we obtain that $\tau$ and $\tau'$ 
are constants and hence $\nu$ is a constant.

{\bf Theorem 3.} Let $M$ be a connected nearly K\"ahler manifold of pointwise constant
antiholomorphic sectional curvature $\nu$ and dim\,$M=2n>4$. Then $M$ is locally 
isometric to one of the following manifolds:

1) The complex Euclidean space $\C^n$;

2) The complex projective space $\C\mathbb P^n$;

3) The complex hyperbolic space $\C\D^n$;

4) The six sphere $\Sph^6$.

{\bf Proof.} From (5) and (7) it follows that 
$$
	2(\nabla_XS)(Y,Z)=S((\nabla_XJ)Y,JZ)+S(JY,(\nabla_XJ)Z) \ ;  \leqno (10)
$$
$$
	(\nabla_XS)(Y,Z)+(\nabla_{JX}S)(JY,Z)=0 \ ;  \leqno (11)
$$
$$
	(\nabla_XS)(Y,Z)+(\nabla_YS)(Z,X)+(\nabla_ZS)(X,Y)=0 \ .  \leqno (12)
$$

Let $\{ E_1,\hdots,E_{2n} \}$ be an arbitrary orthonormal frame field. Taking
into account (10) and (11), from (6) and Theorem 2 we obtain
$$
	\sum_{i=1}^{2n}(\nabla_{E_i}R)(X,Y,Z,E_i)=-\frac16\{(\nabla_XS)(Y,Z)-(\nabla_YS)(X,Z)\}.
$$
This equality and (1) imply $(\nabla_XS)(Y,Z)=(\nabla_YS)(X,Z)$. Then (12) gives that 
the Ricci tensor is parallel: $\nabla_XS=0$.

If $M$ is irreducible, then it is an Einstein manifold and from (9) it follows
that $M$ is of pointwise constant holomorphic sectional curvature. Now the assertion 
follows from the classification theorem in \cite{G1}.

Let $M$ be reducible and be locally a product $M_1(\lambda_1)\times\hdots\times M_k(\lambda_k)$
$(k\ge2)$ where $S=\lambda_ig$ on $M_i(\lambda_i)$ and $\lambda_i$ are different constants
$(i=1,\hdots,k)$. All $M_i$ are nearly K\"ahler manifolds \cite{G2}. Let $X$ on $M_i$,
$Y$ on $M_j$ $(i\ne j)$ be unit vector fields. From $R(X,Y,Y,X)=0$ we find $\nu=0$.
From $R(X,JX,JY,Y)=0$ and (6) it follows that 
$$
	S(X,X)+S(Y,Y)=0 \ .    \leqno (13)
$$ 
If $k>2$, (13) implies $\lambda_i=0$ for $i=1,\hdots,k$. Hence $M$ is of zero Ricci curvature.
The equality (9) gives that $M$ is of zero holomorphic sectional curvature. According to
\cite{G1}, $M$ is locally isometric to $\C^n$. Let $k=2$ and dim\,$M_1=2n_1\ge4$. From
(9) and $\nu=0$ it follows that $M_1$ is of constant holomorphic sectional curvature. Hence
$M_1$ is locally isometric to $\C^{n_1}$. Now (13) gives that $M$ is of zero Ricci curvature
and consequently $M$ is locally isometric to $\C^n$.

\vspace {0.2cm}
\begin{flushright}
{\it Institute of Mathematics  \\
Bulgarian Academy of Sciences \\
Sofia, Bulgaria}
\end{flushright}

\end{document}